# Parallel contact-aware simulations of deformable particles in 3D Stokes flow

Libin Lu[a], Abtin Rahimian[b], Denis Zorin[a]

[a]*Courant Institute of Mathematical Sciences, New York University, New York, NY 10003*
[b]*Department of Computer Science, University of Colorado, Boulder, CO 80302*

**Abstract**

We present a parallel-scalable method for simulating non-dilute suspensions of deformable particles immersed in Stokesian fluid in three dimensions. A critical component in these simulations is robust and accurate collision handling. This work complements our previous work [L. Lu, A. Rahimian, and D. Zorin. Contact-aware simulations of particulate Stokesian suspensions. Journal of Computational Physics 347C: 160–182] by extending it to 3D and by introducing new parallel algorithms for collision detection and handling. We use a well-established boundary integral formulation with spectral Galerkin method to solve the fluid flow.

The key idea is to ensure an interference-free particle configuration by introducing explicit contact constraints into the system. While such constraints are typically unnecessary in the formulation they make it possible to eliminate catastrophic loss of accuracy in the discretized problem by preventing contact explicitly. The incorporation of contact constraints results in a significant increase in stable time-step size for locally-implicit time-stepping and a reduction in the necessary number of discretization points for stability.

Our method maintains the accuracy of previous methods at a significantly lower cost for dense suspensions and the time step size is independent from the volume fraction. Our method permits simulations with high volume fractions; we report results with up to 60% volume fraction. We demonstrated the parallel scaling of the algorithms on up to 16K CPU cores.

## 1. Introduction

Suspensions of rigid and deformable particles are ubiquitous in biological systems and industrial applications. Well-known examples of such fluids are colloids, particulate suspensions, soft-particle pastes, cytoplasm, and blood. The rheology of some these suspensions is well understood in the dilute regime where theoretical models are tractable and simplifying assumptions can be applied to both theoretical and computational models.

In these suspensions, as the volume fraction of particles increases, collective motion emerges and the length scale of the problem grows [24], whereby making theoretical analysis intractable and high-fidelity computational models fraught with numerical challenges. A critical component in simulating *non-dilute* suspensions of (deformable) particles is robust and accurate collision handling. From the rheological perspective, proper collision handling is key for obtaining accurate bulk properties as well as observing correct phase transitions (e.g., to clogging). Computationally, robust collision handling algorithms are key for stable long-time simulations of large ensembles.

To this end, we present an efficient, accurate, and robust method for parallel simulation of dense suspensions in Stokesian fluid in 3D (e.g., Fig. 1). This work complements our previous work [28] by extending it to 3D. In contrast to the 2D case, required computations in 3D are significantly more expensive, and parallelization is essential for for modeling large numbers of particles. In this paper, we introduce new *parallel* algorithms for collision detection and handling in the context of solving integro-differential equations for vesicle flows.

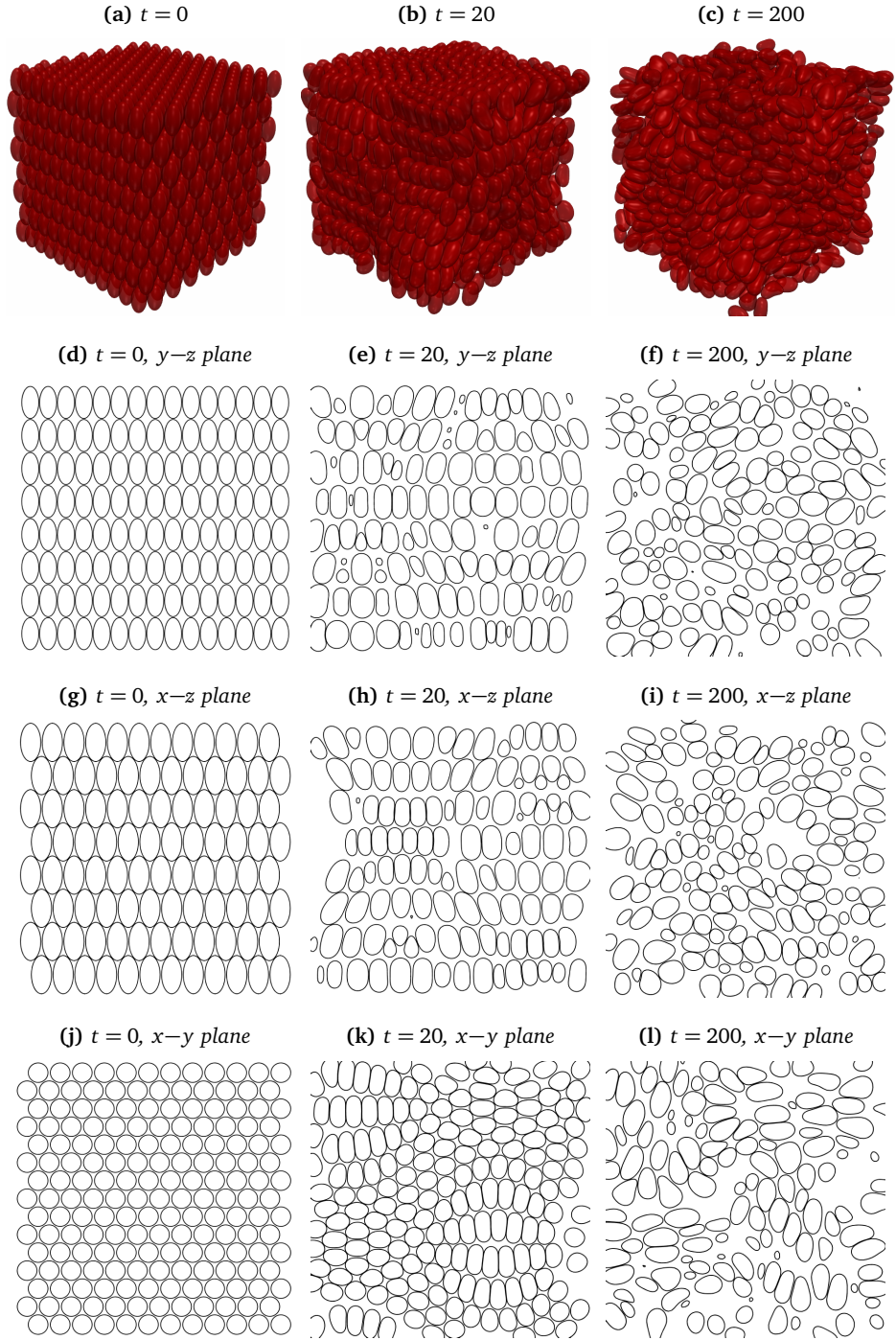

**Figure 1:** TAYLOR VORTEX FLOW EXPERIMENT. *Snapshots of 1440 vesicles in Taylor vortex flow. The volume fraction of vesicles is 58% and each vesicle has a reduced volume of 0.91. For this simulation, we use $p = 16$ spherical harmonic discretization and we upsample the grid to $p = 32$ for contact resolution. We use time step size $\Delta t = 0.1$ and tolerance $1e-5$ for the GMRES solver. (d) to (l) show slice cuts of the flow by different planes. The slice cuts at $t = 0$ look denser because they go through the vesicles' axes of symmetry. The average number of contacts per vesicle per time step is about 2 in this simulation.*

We report computational experiments with volume fractions up to 60% (above the volume fraction of red blood cells, i.e., 45%) involving thousands of deformable particles. We also report strong and weak scaling results for the parallel algorithms on the Stampede system at Texas Advanced Computing Center.

Similar to [28], we focus on the Stokes flow and the *vesicle* model for deformable particles. Vesicles are closed deformable membranes suspended in a viscous fluid that resist bending. They are used to understand the properties of biomembranes [14, 33] and to simulate the motion of blood cells, in which vesicles with moderate viscosity contrast are used to model red blood cells and high viscosity contrast vesicles or rigid particles are used to model white blood cells [3].

While our implementation uses the vesicle deformation model, our contact algorithms do not make any assumptions about the nature of the particles and are applicable to deformable particles with any constitutive properties.

Our approach in this paper and previous work is based on the boundary integral formulation of the Stokes equations. It offers a natural approach for accurate simulation of vesicle flows by reducing the problem to solving equations on surfaces and eliminating the need for discretizing rapidly evolving 3D fluid volumes. However, in non-dilute suspensions, these methods are hindered by a number of difficulties, such as inaccuracies in computing near-singular integrals and artificial force singularities caused by non-physical intersection of particles. In the case of high volume-fraction suspensions, numerical difficulties related to contact and near-contact are of particular importance.

In this work, we resolve contacts explicitly, following [28], we augment the governing equations with a collision-free constraint, i.e., Eq. (2.5).

We briefly recap the reasons for this choice. In principle, in an accurately resolved flow, fluid forces prevent the contact between particles and no explicit collision handling is needed. However, solely relying on the hydrodynamics to prevent contact requires the accurate solution of the flow in the lubrication film, which in turn entails extremely fine spatial and temporal resolution accompanied by increasingly ill-conditioned linear systems in the boundary integral setting [47, 49] — becoming computationally impractical as the volume fraction increases. In comparison, explicit contact handling can ensure robustness and accuracy of simulations, without imposing excessively restricting constraints on the spacial discretization and time stepping.

Contact is often resolved by introducing artificial repulsion forces along with adaptive time stepping [31, 45, 46] This approach works well for dilute suspensions. However, the time-step for these methods is determined by the closest pair of vesicles, and tends to be uniformly small for dense suspensions (see Section 4.5). Furthermore, the heuristic parameters of the repulsion force may be difficult to determine automatically, and often require problem-dependent adjustments (Section 4.5 and Table 5).

Compared to repulsion, contact constraints ensure that the geometry remains intersection-free, even for relatively coarse spatial and temporal discretizations, where the fidelity of the numerical model is insufficient for resolving the lubrication film precisely enough to prevent contact, without any simulation-specific parameter tuning.

*Our contributions.* We generalize the method of [28] to three-dimensional flows, specifically, changing the formulation to use four-dimensional space-time volumes for collision detection, with mesh, instead of piecewise-linear proxies. One of our principal contributions is the distributed parallel version of contact detection and resolution algorithms. Our algorithm for contact detection, based on fast parallel sorting [55], utilizes an implicit grid, sorted in Morton curve order on distributed memory machines, for inter-process collision detection, and uses an explicit spatial grid to perform intra-process collision detection calculation locally. The contact constraints are formulated as a Nonlinear Complementarity Problem (NCP). The solution of the NCP problem and the construction of the necessary matrices are both done in parallel based on the iterative minimal map Newton method.

We accelerate the computation of far-field hydrodynamics interactions using the highly optimized parallel PVFMM library [29]. The PVFMM library supports periodic boundary conditions and this allows us to simulate vesicles in periodic flows with a prescribed volume fraction.

Our method makes possible simulations with high volume fractions (we report results with up to 60% volume fraction, e.g., Fig. 1). Within our framework, the time step size is independent of the volume fraction and the simulation wall-clock-time is at least an order of magnitude faster than the adaptive time step approach (Section 4.5).

*Related work.* This paper is most closely related to [28, 31, 49]. We extend these works to enable long-time contact-aware simulation of concentrated vesicle suspensions in parallel. In [28], we proposed the initial version of contact detection and resolution algorithms for 2D and demonstrated that these algorithms enable long-term simulations, increase robustness, and reduce the computational costs. In [31], several novel computational algorithms (e.g., adaptive time stepping, robust near-singular integration, and adaptive mesh refinement) are proposed to facilitate long-time simulation of concentrated suspensions. A short-range repulsion force along with adaptive time stepping were introduced to avoid collision, permitting simulations with up to 35% volume fraction. If a collision is detected, the time step is backtracked and time step size is refined, which in high volume-fraction cases, may result in very small time steps (Section 4.5). Our method enables higher density simulations with significantly larger time steps, and does not require parameter tuning for the penalty function.

More broadly, Stokesian particle models are employed to theoretically and computationally investigate the properties of biological membranes [51], drug-carrying capsules [54], and blood cells [34, 40]. There is an extensive body of work on numerical methods for Stokesian particulate flows and a review of the literature up to 2001 can be found in [42]. Reviews of later advances can be found in [48, 49, 60]. Here, we briefly summarize some notable numerical methods and discuss the most recent developments.

Integral equation methods have been used extensively for the simulation of Stokesian particulate flows such as droplets and bubbles [26, 27, 50, 70], vesicles [10, 13, 40, 48, 49, 53, 60, 67, 68], and rigid particles [38, 39, 66]. Other methods — such as phase-field approach [4, 7], immersed boundary and front tracking methods [21, 65], and level set method [23] — are used by several authors for the simulation of particulate flows.

In [28] we extensively reviewed works on collision detection and handling in the context of (i) Stokesian flows [13, 26, 35, 46, 52, 69, 71, 72]; (ii) contact mechanics and response [12, 19, 22, 44, 59, 63, 64]; and (iii) computer graphics [2, 8, 11, 15, 16, 17, 36, 43] . We refer the interested reader to [28] for more details; here, we focus on most closely related work.

Our constraint-handling algorithms belongs to a large family of *constraint-based methods*, commonly used to handle contacts reliably in many applications, primarily in computer graphics. This class of methods meets our goals of providing robustness and improving efficiency of contact response, while minimizing the impact on the physics of the system. Our contact resolution approach is directly based on [17] and is closest to [1], in which the intersection volume and its gradient with respect to control vertices are computed at the candidate step. Harmon et al. [17] assumes linear trajectory between edits and defines a *space-time interference volume* (STIV) which serves as a gap function. We refine this formulation to define our contact constraint.

Parallel algorithms for collision detection, for shared and distributed memory architectures, often include two stages: the first stage is culling, i.e., reducing the set of potential collision pairs using a fast and conservative criterion to determine which pairs do no intersect, followed by precise collision detection between remaining pairs. Many authors used multi-threading and GPU computation to accelerate different stages of collision detection [25, 32, 56, 57, 58]. The focus of these algorithms is distributing the task of collision culling between threads, for which they use a variety of techniques.

[25] proposed a GPU-based discrete collision detection algorithm, in which axis-aligned bounding boxes are computed for each object followed by sweep-and-prune steps on the GPU to identify a small set of collision candidates efficiently. In Mazhar, Heyn, and Negrut [32], surfaces are approximated by large collections of padded spheres and intersections between spheres are used to cull collision candidates. Kim et al. [20]

| Symbol | Definition | Symbol | Definition |
|---|---|---|---|
| $\gamma_i$ | The boundary of the $i^{\text{th}}$ vesicle | LI | Locally Implicit time stepping |
| $\gamma$ | $\cup_i \gamma_i$ | CLI | *Constrained* Locally Implicit time stepping |
| $\mu$ | Viscosity of the ambient fluid | GI | Globally Implicit time stepping |
| $\mu_i$ | Viscosity of the fluid inside $i^{\text{th}}$ vesicle | RGI | Repulsion-based Globally Implicit time stepping |
| $\nu_i$ | The viscosity contrast $\mu_i/\mu$ | | |
| $\sigma$ | Tension | $d_m$ | Minimum separation distance |
| $\chi$ | Shear rate | $\boldsymbol{f}_\sigma$ | Tensile force |
| $\varpi_i$ | The domain enclosed by $\gamma_i$ | $\boldsymbol{f}_b$ | Bending force |
| $\varpi$ | $\cup_i \varpi_i$ | $\boldsymbol{f}_c$ | Collision force |
| $\mathscr{G}$ | Stokes Single-layer operator | $J$ | Jacobian of contact volumes $V$ |
| $\mathscr{T}$ | Stokes Double-layer operator | $\boldsymbol{n}$ | Unit outward normal |
| LCP | Linear Complementarity problem | $\boldsymbol{u}$ | Velocity |
| NCP | Nonlinear Complementarity Problem | $\boldsymbol{u}^\infty$ | The background velocity field |
| STIV | Space-Time Interference Volumes | $V$ | Contact volumes |
| | | $\boldsymbol{X}$ | A Lagrangian point on a surface |

**Table 1:** INDEX OF FREQUENTLY USED SYMBOLS, OPERATORS, AND ABBREVIATIONS.

developed a parallel continuous collision detection algorithm for heterogeneous shared-memory architectures. The algorithm uses bounding volume hierarchy (BVH) for culling. CPUs perform the BVH traversal and culling, while GPUs perform collision tests between geometric primitives.

For distributed memory collision handling, more complex data structures are required. Iglberger and Rüde [18] presents a method for rigid body simulation on distributed memory machines. The domain of interest is partitioned and distributed over all MPI processes and each process maintains a list of objects in its domain.

The algorithms we present in this work shares a number of features with [6, 37, 62]. Warren and Salmon [62] presents a parallel distributed-memory N-Body algorithm. This algorithm uses a parallel octree using the Morton ordering curve and a local essential tree to form a distributed memory octree for N-body interaction calculation. Du et al. [6] describes a parallel continuous collision detection algorithm for rigid bodies. Rigid bodies are approximated by spherical bounding volumes and space-time axis aligned bounding boxes are computed for each sphere. The domain is divided into cells, each with a list of overlapping rigid bodies. The cells are distributed over MPI processes and the lists are dynamically updated in parallel. Within each cell, the collision detection is performed locally. Similarly, Pabst, Koch, and Strasser [37] (a shared-memory algorithm) use spatial cells along with hashing for candidate collision identification, followed by primitive collision tests.

*1.1. Nomenclature*

In Table 1, we list symbols and operators used in this paper. Lowercase letters refer to scalars, and lowercase bold letters refer to vectors. Discretized quantities are denoted by sans serif letters.

*1.2. Synopsis of the method*

We use the boundary integral formulation based on [49]. The basic formulation uses integral equation form of the problem and includes the effects of the viscosity contrast. We add contact constraints to this formulation as an inequality constraint on a gap function that is based on *space-time intersection volume* [28], the 3D *space-time intersection volume* is presented in Appendix A. The contact force is then parallel to the gradient of this volume with the Lagrange multiplier as its magnitude.

We solve the resulting contact NCP for the Lagrange multipliers of the constraints using a Newton-like matrix-free method as a sequence of Linear Complementarity Problems (LCP) [5, 9]. Each LCP is solved iteratively using GMRES. The spherical harmonics bases are used for spatial discretization. For time stepping, we use semi-implicit backward Euler (Section 3).

In Section 3, we present a summary of the spatial and temporal discretization. We also present parallel algorithms for collision culling, STIV computation, and the solution of LCP's.

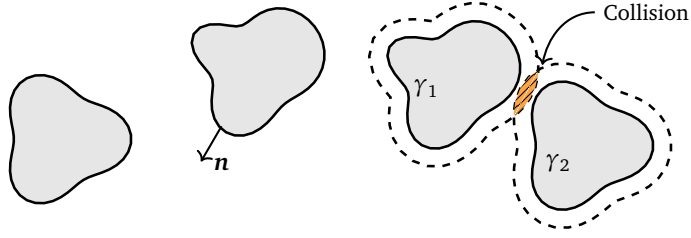

**Figure 2:** SCHEMATIC. *Vesicles are suspended in free-space, both filled with fluid. The vesicle boundaries are denoted by $\gamma_i$ ($i = 1, \ldots, N$). $\mathbf{n}$ denotes the outward normal vector to the vesicles' boundaries. The dotted lines around vesicles' boundaries denote the prescribed minimum separation distance for each of the vesicle. The minimum separation distance is a parameter can be set to zero. In this schematic, vesicles $\gamma_1$ and $\gamma_2$ are in contact. A space-time intersection volume is marked by ▨.*

In Section 4, we present results showing the accuracy and effectiveness of our scheme. We also present weak and strong scaling results for simulations on up to 16K processors and report wall-clock-time for different methods and volume fractions.

## 2. Formulation

In this section, we briefly outline the mathematical formulation of vesicle suspensions in Stokes flow with contact constraints. For a more in-depth formulation applicable to 2D and 3D vesicle flows, see [28, 49].

We consider $N$ vesicles suspended in an incompressible Stokes fluid. The fluid domain is assumed to be unbounded. The governing equations have the form

$$-\mu \Delta \mathbf{u}(\mathbf{x}) + \nabla p(\mathbf{x}) = \mathbf{F}(\mathbf{x}) \quad \text{and} \quad \nabla \cdot \mathbf{u}(\mathbf{x}) = 0 \qquad (\mathbf{x} \in \mathbb{R}^3), \tag{2.1}$$

$$\mathbf{F}(\mathbf{x}) = \int_\gamma \mathbf{f}(\mathbf{X}) \delta(\mathbf{x} - \mathbf{X}) \, dA(\mathbf{X}), \tag{2.2}$$

where $\mu$ denotes the viscosity of ambient fluid and $\mathbf{f}$ is the surface density of the force exerted by the vesicles' membrane. We let $\mathbf{x}$ denote an Eulerian point in the fluid ($\mathbf{x} \in \mathbb{R}^3$) and $\mathbf{X}$ denotes a Lagrangian point on a vesicle. We assume the no-slip boundary condition on the surface of vesicles

$$\mathbf{u}(\mathbf{X}, t) = \mathbf{X}_t \qquad (\mathbf{X} \in \gamma). \tag{2.3}$$

The vesicle membrane is assumed to be inextensible, i.e.,

$$\nabla_\gamma \cdot \mathbf{u}(\mathbf{X}) = 0 \qquad (\mathbf{X} \in \gamma), \tag{2.4}$$

where $\nabla_\gamma \cdot$ denotes the surface divergence operator.

### 2.1. Contact constraint

We follow the approach of [28] to address the contact issue resulting from inadequate computational accuracy. We enforce contact constraint in the system, i.e.,

$$V(\gamma, t) \geq 0. \tag{2.5}$$

We derive the 3D contact constraint function (space-time interference volume) in Appendix A. To resolve the contact constraint in parallel in a distributed memory architecture is one of the main challenges in practical application of this approach. In Section 3.3, we describe in detail the parallelization of contact constraint computation and other parallizations needed. The equations of motion (Eqs. (2.1), (2.3), and (2.4)), along with the constraint Eq. (2.5) yield the modified equations:

$$-\mu \Delta \boldsymbol{u} + \nabla p = \boldsymbol{F}', \tag{2.6}$$

$$\boldsymbol{F}'(\boldsymbol{x}) = \boldsymbol{F}(\boldsymbol{x}) + \int_\gamma (\boldsymbol{f}_b + \boldsymbol{f}_\sigma + \boldsymbol{f}_c)\delta(\boldsymbol{x} - \boldsymbol{X})\,\mathrm{d}A, \tag{2.7}$$

$$\boldsymbol{f}_b = -\kappa_b \left[\Delta_\gamma H + 2H(H^2 - K)\right]\boldsymbol{n} \tag{2.8}$$

$$\boldsymbol{f}_\sigma = \sigma \Delta_\gamma \boldsymbol{X} + \nabla_\gamma \sigma, \tag{2.9}$$

$$\boldsymbol{f}_c = (\mathrm{d}_{\boldsymbol{u}} V)^T \lambda, \tag{2.10}$$

$$\lambda \geq 0; \lambda \cdot V = 0, \tag{2.11}$$

where $\mathrm{d}_{\boldsymbol{u}} V$ is the variation of $V$ with respect to $\boldsymbol{u}$, $\kappa_b$ is the membrane's bending modulus, $\boldsymbol{n}$ is the outward normal (as shown in Fig. 2) to the vesicle surface $\gamma_i$, $H$ and $K$ are respectively mean and Gaussian curvatures, and the collision force $\boldsymbol{f}_c$ is added to the traction jump across the vesicle's interface.

It is conventional to write $V \geq 0$, $\lambda \geq 0$, and $\lambda \cdot V = 0$, into one complementarity condition expression

$$0 \leq V(t) \quad \perp \quad \lambda \geq 0. \tag{2.12}$$

The formulae outlined above govern the evolution of the vesicle. Following the same approach in [28], we use boundary integral formulation to solve the problem. For the sake of completeness, we briefly summarize the boundary integral formulation in Appendix B. Given the configuration of the vesicle, the unknowns are velocity $\boldsymbol{u}(\boldsymbol{X})$, tension $\sigma$ of vesicles' interface determined by Eqs. (B.5–B.7) and $\lambda$ which results from the Signorini (KKT) conditions for the contact constraint, Eq. (2.12) ($\lambda$ determines the strength of the contact force). The velocity is integrated for the vesicles' trajectory using no-slip boundary condition Eq. (2.3).

## 3. Numerical algorithms

In this section, we describe the algorithms required for solving the dynamics of particulate Stokesian suspensions. Our spatial/temporal representation and surface quadratures follow our previous work [30, 49, 61].

We present a set of parallel distributed-memory algorithms for contact resolution in 3D. The parallelization is a necessity for simulations of reasonable size in 3D. The overall approach follows our two-dimensional contact-aware scheme presented in [28], but new algorithms are needed for the distributed-memory contact detection and resolution.

At every time step, we resolve contacts by solving a nonlinear complementarity problem (NCP) in parallel. The NCP is solved iteratively by recursive linearization and using a new parallel Linear Complementarity Problem (LCP) solver.

In the following sections, we first summarize a brief description of the spatial discretization from [49], then discuss the time discretization with contact constraint and the corresponding parallel algorithms.

### 3.1. Spatial discretization

*3.1.1. Surface representation and Quadrature rule.* The boundary of each vesicle is assumed to be parameterized by a smooth map $\boldsymbol{X}(\phi, \theta)$ from a unit sphere $\mathbb{S}^2$ to $\mathbb{R}^3$. We use spherical harmonics expansions to represent the surface and all functions on the surface of vesicle (such as surface position/velocity and tensile/bending/contact forces) [49, 61].

We use the same quadrature rule to evaluate single- and double-layer potential in our boundary integral formulation. For the target point $\boldsymbol{Y}$ not on the surface $\gamma$, the integrand is smooth; we use Gauss-Legendre quadrature for $\phi$ and trapezoidal quadrature for $\theta$ directions. We use the Fast Multipole Method (FMM) to accelerate the calculation. As $\boldsymbol{Y}$ approaches $\gamma$, the quadrature rule is insufficient to accurately evaluate near-singular layer-potential. We use the near-singular integration scheme discussed in [30]. For the target point $\boldsymbol{Y}$ on the surface $\gamma$, the layer potential is singular. We use the singular integration algorithm discussed in [30, 61] for both the Stokes single-layer potential and the Stokes double-layer potential. This quadrature scheme involves rotation of the spherical harmonics and is spectrally accurate for both single- and double-layer potentials.

*3.1.2. Piecewise-linear triangular discretization for constraints.* Similar to the 2D case, while the spectral spatial discretization is used for most computations, it poses a problem for the minimal-separation constraint discretization. Computing parametric surface intersections, an essential step in the STIV computation, is relatively expensive and difficult to implement robustly, as this requires solving nonlinear equations for intersections. We use a piecewise-linear triangular discretization of the surface to calculate the space-time contact volume and its gradient.

The collocation points for the $q$-grid naturally give a quadrilateral mesh of the surface, except at the poles, where quads are degenerate. Splitting each quadrilateral along the diagonal we get a piecewise-linear triangular discretization of the surface as illustrated in Fig. 3.

We first upsample the $q$-grid to a $q^{up}$-grid and convert the $q^{up}$-grid to a triangular mesh for STIV calculation. Since we opt for a low-order, piecewise-linear approximation of the vesicle, we need to use an algorithm that ensures that at least the target minimal separation is maintained between actual smooth surfaces represented by spherical harmonics. For the triangular mesh, we set the separation distance to $(1+2\alpha)d_m$, where $d_m$ is the target minimum separation distance. We observe that the sensitivity to the separation distance on the overall accuracy is low in most situations as explored in Section 4; in practice we choose $\alpha = 0.05$ and we find that $q^{up} = 32$ is sufficient for $\|\mathbf{X}^{q^{up}} - X(\phi,\theta)\|_\infty < \alpha d_m$.

To compute the contact constraint with the triangular mesh of a $q^{up}$-grid, we calculate the discretized space-time contact volume as the sum of triangle-vertex contact volumes $V = \sum_k V_k(\mathbf{t}, \mathbf{X})$, where $k$ indexes triangle-vertex pairs. Parallel algorithms to find colliding triangle-vertex pairs are discussed Section 3.3.

For each triangle-vertex pair $(\mathbf{t}(\mathbf{X}_{i-1}, \mathbf{X}_i, \mathbf{X}_{i+1}), \mathbf{X}_k)$ pair, we solve a degree six equation to find their earliest contact time $\tau_I$ assuming linear trajectory between the initial position at time $t_n$ and candidate position at $t_{n+1}$:

$$([\mathbf{X}_k(t) - \mathbf{X}_{i-1}(t)] \cdot [(\mathbf{X}_i(t) - \mathbf{X}_{i-1}(t)) \times (\mathbf{X}_{i+1}(t) - \mathbf{X}_{i-1}(t))])^2 \\ - d_m^2 \left\| [(\mathbf{X}_i(t) - \mathbf{X}_{i-1}(t)) \times (\mathbf{X}_{i+1}(t) - \mathbf{X}_{i-1}(t))] \right\|^2 = 0, \quad (3.1)$$

where $\mathbf{X}_k(t) = \mathbf{X}_k(t_n) + t\mathbf{U}_k$, with the linear trajectory assumption the vertex velocity is defined as $\mathbf{U}_k = [\mathbf{X}_k(t_{n+1}) - \mathbf{X}_k(t_n)]/\Delta t$. We calculate the triangle-vertex contact volume as follows (Appendix A):

$$V_k(\mathbf{t}, \mathbf{X}) = (t - \tau_I)(\epsilon^2 + (\mathbf{U}_k \cdot \mathbf{n}(\tau_I))^2)^{1/2}|\mathbf{t}|, \quad (3.2)$$

where $\mathbf{n}(\tau_I)$ is the normal to the triangle $\mathbf{t}$ at the intersection time and $|\mathbf{t}|$ is the area of the triangle. For each triangle-vertex contact volume, we calculate the gradient of Eq. (3.2) with respect to the vertices $\mathbf{X}_{i-1}$, $\mathbf{X}_i$, $\mathbf{X}_{i+1}$ and $\mathbf{X}_k$ and by summing over all the triangle-vertex contact pairs, we get the total space-time interference volume and its gradient.

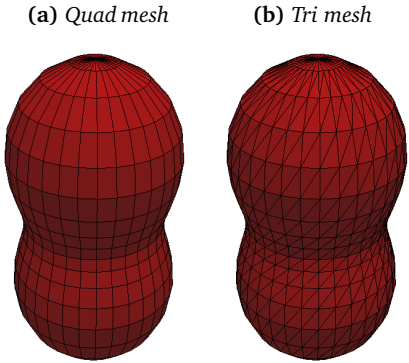

**Figure 3:** PIECEWISE-LINEAR TRIANGULAR DISCRETIZATION. *Illustrating the conversion of a 16-grid (a) to a piecewise-linear triangular mesh (b).*

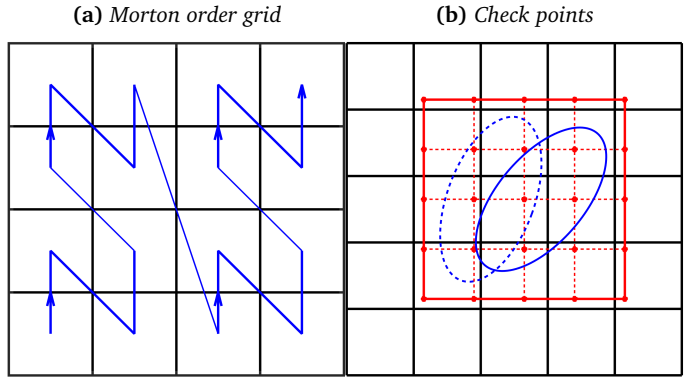

**Figure 4:** MORTON ORDER GRID. *Diagram (a) shows the 2D Morton curve order for a grid. (b) shows the schematic of 2D space-time bounding box (red) of a vesicle (dotted blue curve is the initial position, solid blue curve is the candidate position) filled with check points.*

*3.2. Temporal discretization*

Our temporal discretization is based on the locally-implicit time-stepping scheme in [49]: we treat intra-vesicle interactions implicitly and inter-vesicle interactions explicitly. Following [28], we combine the first-order backward Euler time stepping with minimal-separation constraint. More accurate time-stepping methods such as spectral deferred correction (SDC) method [28] can be easily integrated. We denote this locally-implicit scheme with collision constraint by CLI. In Section 4, we compare this method to the same scheme without constraints (LI) and the globally semi-implicit with/without repulsion(GI and RGI) schemes, where all interactions treated implicitly as in [30].

Marking the unknowns to be solved for with '·+' superscript, we discretize Equations (2.6–2.11) and obtain the time-stepping equation in a compact form as

$$\mathbf{A}\mathbf{X}^+ = \mathbf{b} + \mathbf{G}\mathbf{f}_c^+, \tag{3.3}$$

$$0 \leq \mathsf{V}(\gamma; t^+) \quad \perp \quad \lambda \geq 0. \tag{3.4}$$

$\mathbf{A}$ is a block diagonal matrix, each block $\mathbf{A}_{ii}$ of $\mathbf{A}$ corresponds to the self-interaction of $i^{\text{th}}$ vesicle. $\mathbf{G}$ is the discretized Stokes operator, it is also a block diagonal matrix. The right-hand side $\mathbf{b}$ includes all the explicitly treadted inter-vesicle interactions. At each time step, we resolve contacts by solving this mixed Nonlinear Complementarity Problem (NCP). We outline our parallel algorithms for solving this NCP in the next section.

To solve the NCP Eqs. (3.3) and (3.4), we linearize the nonlinear constraint function $\mathsf{V}(\gamma; t)$ around the current candidate position and iteratively solve a sequence of LCPs:

$$\mathbf{A}\triangle\mathbf{X}^{+,k} = \mathbf{G}\mathbf{J}^T \lambda^{+,k}, \tag{3.5}$$

$$0 \leq \mathsf{V}(\gamma; t^{+,k}) + \mathbf{J}\triangle\mathbf{X}^{+,k} \quad \perp \quad \lambda^{+,k} \geq 0, \tag{3.6}$$

until the original NCP is solved to the desired accuracy. In Eq. (3.5) and Eq. (3.6), $\triangle\mathbf{X}^{+,k}$ is the update to get the new candidate solution in $k^{\text{th}}$ contact-resolving iteration, and $\mathbf{J}$ denotes the Jacobian of the STIV $d_X \mathsf{V}(\gamma, t^{+,k})$, where $k$ indexes the contact-resolving iterations. We call each LCP solve a contact-resolving iteration. After each contact-resolving iteration, the current candidate position is updated with $\triangle\mathbf{X}^{+,k}$ and $\mathbf{J}^T \lambda^{+,k}$ is accumulated into the contact force $\mathbf{f}_c^+$. By using the factor that $\triangle\mathbf{X}^{+,k} = \mathbf{A}^{-1} \mathbf{G}\mathbf{J}^T \lambda^{+,k}$, we can cast Eq. (3.5) and Eq. (3.6) into the standard LCP form solving for $\lambda^{+,k}$

$$0 \leq \mathsf{V} + \mathbf{B}\lambda^{+,k} \quad \perp \quad \lambda^{+,k} \geq 0, \tag{3.7}$$

where $\mathbf{B} = \mathbf{J}\mathbf{A}^{-1}\mathbf{G}\mathbf{J}^T$ and $\mathsf{V} = \mathsf{V}(\gamma, t^{+,k})$ which are calculated using current candidate position.

The top-level iteration described above constitutes our contact-free time-stepping algorithm (cf. [28, Algorithm 1]). All steps of this algorithm are parallelized. The inter-vesicle interactions (explicit part of the time stepping) are computed with a parallel Fast Multipole Method [29]. The STIVs are computed using a parallel collision detection algorithm Section 3.3. We explicitly form the distributed LCP matrix and solve the LCP in parallel, as discussed in Section 3.3.3.

*3.3. Parallel collision handling*

In this section, we describe the most challenging algorithmic part of our method, parallel collision handling, which is essential for scalability.

To avoid costly communication and computation, our contact detection is performed in two phases. In the first phase, we find intersecting bounding boxes of particles (Alg. 1). For each particle, this results in a list of other particles that it may be colliding with (i.e., *candidate pairs*). In the second phase, we communicate the mesh information for the particles in the list and compute the pairwise STIVs (Alg. 2). These intersection volumes and their gradients are used in the LCP to find the magnitude of contact force.

In our algorithm descriptions, superscript $p$ denotes data that resides on process $p$. Each vesicle is assigned to a process $p$ and the assignment does not change during the simulation. Variables without superscripts are either shared variables among all processes, or global arrays whose local parts are denoted with the superscripts. We use subscripts for indices of the vesicles, bounding boxes, etc. Table 2 summarizes the variables used in this section.

| Symbol | Definition |
|---|---|
| $n_p$ | Number of MPI processes |
| $p$ | Process index |
| $n$ | Number of points on a vesicle (assumed fixed) |
| $N^p$ | Number of vesicles |
| $\overline{N}^p$ | Number of ghost vesicles |
| $N_c^p$ | Number of check points |
| $N_P^p$ | Number of candidate vesicle pairs |
| $N_v^p$ | Number of contact volumes |
| $X_0^p, X_1^p$ | Initial and candidate positions (size $nN^p$) |
| $\overline{X}_0^p, \overline{X}_1^p$ | Initial and candidate positions of ghost vesicles |
| $I^p$ | Set of global indices for vesicles (vectors of size $N^p$) |
| $\overline{I}^p$ | Set of global indices for ghost vesicles |
| $B^p$ | Set of space-time bounding boxes of vesicles |
| $I_b^p$ | Set of global bounding-boxes indices |
| $P^p$ | Set of index pairs of intersecting bounding boxes, $\{(j_r, k_r) \mid r = 1, \ldots, N_P^p, j_r \in I_b^p, k_r \in I_b, j_r < k_r\}$ |
| $\mathbf{c}^p, I_c^p, I_{c,b}^p, M^p$ | Positions, global indices, vesicles indices, and Morton codes of check points |
| $h_m$ | Morton order grid cell size |
| $V^p, I_v^p$ | Contact volumes and their global indices (size $N_v^p$) |
| $J^p$ | Jacobian of the vector of contact volumes with respect to vesicles (size $N_v^p \times 3nN^p$) |
| $\overline{J}^p$ | Jacobian with respect to ghost vesicles (size $N_v^p \times 3n\overline{N}^p$) |
| $\lambda^p = \lambda(I_v^p)$ | Lagrange multipliers |

**Table 2:** PARALLEL VARIABLES. *Superscript p denotes the data corresponding to process p.*

*3.3.1. Phase 1: Bounding box intersections.* To narrow down the set of potential collision pairs efficiently, we initially use *space-time* bounding boxes as collision proxies, i.e., 3D bounding boxes enclosing point trajectories from the initial positions $X_0^p$ to candidate positions $X_1^p$ for each vesicle. Each axis-aligned bounding box is stored as a pair of points $\{\underline{\mathbf{b}}_i^p, \bar{\mathbf{b}}_i^p\}$, that are the corners of the box with lexicographically-ordered minimum and maximum coordinate values.

To identify the intersecting bounding box pairs efficiently and in parallel, we use a spatial grid algorithm. There are three main steps in finding the intersecting bounding box pairs using a spatial grid:
(1) For each bounding box, find grid cells it overlaps.
(2) For each grid cell, compute a list of bounding boxes overlapping it by merging the lists of (box, grid cell) pairs from step 1, and sorting it by grid cell.
(3) For each grid cell with a non-empty list, perform intersection check for all bounding box pairs in that cell and find all the intersecting bounding box pairs.

*Parallel version.* The most direct approach to parallelizing this algorithm is to distribute the lists of boxes associated with grid cells across processes. However, performing step (2) efficiently is difficult in this case, since merging the lists obtained in step (1) requires irregular communication pattern between all processors.

Instead of using an explicitly distributed spatial grid data structure storing lists of boxes, we use an implicit representation of grid cells based on Morton curve order numbering. Morton IDs of cells are assigned to *check points*, sampled on bounding boxes in a way that guarantees that at least one check point is contained in every grid cell overlapping the box. Parallel-sorting the check points by their Morton IDs collects, on each processor, a set of check points corresponding to bounding boxes overlapping the same grid cells, as their Morton IDs will be the same. Then bounding boxes overlapping each grid cell are checked for intersections in parallel. Finally, detected candidate intersection pair lists are scattered back to the processors owning bounding boxes contained in the list. Next, we describe the algorithm more formally. In the pseudocode (Alg. 1) lines 1–3, line 4, and lines 5–11 respectively are steps (1) to (3) outlined above.

**Algorithm 1:** BOUNDING BOX INTERSECTION.

    **input**  : A set of axis aligned bounding boxes $B^p$, $I_b^p$
    **output:** Intersecting bounding box index pairs: $P^p$

    // (1) Insert boxes into the tree by computing the Morton ID of grid cells they intersect.

1  $h_m \leftarrow \texttt{average}\left(\|\underline{\mathbf{b}}_i^p - \bar{\mathbf{b}}_i^p\|\right)$                                       // Parallel reduction
2  $\{\mathbf{c}^p, I_c^p, M^p, I_{c,b}^p, \underline{\mathbf{x}}^p, \bar{\mathbf{x}}^p\} \leftarrow \texttt{generateCheckPts}\left(B^p, I_b^p, h_m\right)$
3  Update $I_c^p$ to global index                                                               // Parallel scan

    // (2) Merge the lists by sorting check points using Morton IDs and scatter data based on $I_c^p$.

4  $\texttt{hypercubeQSort}\left(M^p, \{I_c^p, I_{c,b}^p, \underline{\mathbf{x}}^p, \bar{\mathbf{x}}^p\}\right)$

    // (3) Local check of intersecting bounding box pairs

5  $P^p \leftarrow \emptyset$
6  **for** *each* $m \in M^p$ **do**
7     **for** *each bounding box pair* $\{j,k\}(j,k \in I_{c,b}^p)$ *in grid* $m$ **do**
8         **if** $\texttt{bbiCheck}(\underline{\mathbf{b}}_j, \bar{\mathbf{b}}_j, \underline{\mathbf{b}}_k, \bar{\mathbf{b}}_k)$ **then**
9             Add intersecting bounding box pair$\{j,k\}$ to $P^p$

10  Send $P^p(\{j,k\})$ to the process $q$ such that $j \in I_b^q$
11  $P^p \leftarrow \texttt{unique}(P^p)$

|  | Process 1 | | | | Process 2 | | | |
|---|---|---|---|---|---|---|---|---|
|  | Unsorted | | | | | | | |
| $M$: | 111 | 110 | 001 | 000 | 001 | 000 | 000 | 101 |
| $I_c$: | 1 | 2 | 3 | 4 | 5 | 6 | 7 | 8 |
| $I_{c,b}$: | 1 | 1 | 1 | 1 | 2 | 2 | 2 | 2 |
| local data: | a | b | c | d | e | f | g | h |
|  | Sorted (shaded rows are misaligned) | | | | | | | |
| $M$: | 000 | 000 | 000 | 001 | 001 | 101 | 110 | 111 |
| $I_c$: | 4 | 6 | 7 | 3 | 5 | 8 | 2 | 1 |
| $I_{c,b}$: | 1 | 1 | 1 | 1 |  | 2 | 2 | 2 | 2 |
| local data: | a | b | c | d |  | e | f | g | h |
|  | Scattered (data is aligned) | | | | | | | |
| $M$: | 000 | 000 | 000 | 001 | 001 | 101 | 110 | 111 |
| $I_c$: | 4 | 6 | 7 | 3 | 5 | 8 | 2 | 1 |
| $I_{c,b}$: | 1 | 2 | 2 | 1 | 2 | 2 | 1 | 1 |
| local data: | d | f | g | c | e | h | b | a |

**Table 3:** SORTING AND SCATTERING OF CHECK POINTS ON TWO PROCESSES. *A simple check points sorting and data scattering example. This example illustrates the data movement between two processes. Unsorted rows show the input information on process 1 and process 2. Sorted rows show the sorted check points' Morton codes M and shuffled index set $I_c$ with $I_{c,b}$ and unscattered data on each process. Shaded rows are misaligned. Scattered rows show the result of scattering using the shuffled index set $I_c$ as a scatter mapping from original data to scattering place on each process.*

We view the spatial grid as a uniformly refined octree. The depth of octree is $\log(L/h_m)$ where $h_m$ is the grid cell size. We set $h_m$ to the average diagonal length of all vesicle bounding boxes and $L$ is the domain length. Each leaf in the octree is associated with a unique Morton ID, as shown in Fig. 4(a). Table 3 illustrates the data movement in the parallel sorting algorithm(line 4 of Alg. 1).

In step (1), check points $\mathbf{c}^p$ are defined for the bounding boxes $B^p$ (lines 2 and 3). As schematically shown in Fig. 4(b), check points are located uniformly inside each bounding box with a spacing less than $h_m$. With this spacing, a bounding box should have at least one check point in an any grid cell it overlaps. Each check point is associated with its originating bounding box index and that bounding box's spatial data, i.e., $I_{c,b}^p$ and $\{\underline{\mathbf{x}}^p, \bar{\mathbf{x}}^p\}$. Our choice of spacing results in $2^3$ to $3^3$ check points for each bounding box.

We compute the set of all Morton IDs $M^p$ for all check points. Since each Morton ID is uniquely associated with a grid cell and each bounding box has at least one check point inside the grid cell overlapping with that bounding box, our lists $M^p$ include the Morton ID of all grid cells that intersect with each bounding box, possibly multiple times (we drop the duplicates).

Initially, the lists $I_c^p$ are the global indices of the check points on process $p$; to calculate this global index we perform an MPI scan operation on the number of check points on $p^{\text{th}}$ process $N_c^p$. With the Morton IDs and global check point index computed for each check point, we can move to the next stage.

Step (2) of the algorithm (merging lists) is equivalent to a parallel sort on the Morton IDs followed by scattering of data associated with a check point using the shuffle index $I_c^p$ obtained as a result of the sort.

We use a parallel hypercube quicksort algorithm [55]. The hypercube quicksort algorithm takes the Morton code $M^p$ (keys) and index set of check points and their associated information $D_c^p = \{I_c^p, I_{c,b}^p, \underline{\mathbf{x}}^p, \bar{\mathbf{x}}^p\}$. (values) on $p^{\text{th}}$ process as input, it outputs the globally sorted Morton code, and shuffled per-check point information.

The parallel sorting algorithm evenly distributes the sorted Morton IDs among all processes so that the $p^{\text{th}}$ process contains list $M^p$, which is a consecutive part of the sorted global Morton code $M$, along with corresponding data $D_c^p$. In addition, the sorting algorithm places identical Morton codes on the same process, which is equivalent to assigning the list of boxes overlapping a grid cell to a single process. This process is responsible for checking intersections of boxes on this grid cell's list.

We find intersecting bounding box pairs on process $p$ by checking all pairs of bounding boxes in $I_{c,b}^p$ which have the same Morton code (line 5–11 of Alg. 1). Line 8 checks the candidate intersecting bounding box pairs to see whether they actually intersect or not. For each pair of intersecting boxes $(j,k)$ we generate two ordered pairs $(k,j)$ and $(j,k)$.

The process $p$ sends the pairs to the process which originally owned the first vesicle in each $P^p$ entry, using a sparse MPI all-to-all communication call. After the MPI communication, the $p^{\text{th}}$ process will have a list of intersecting bounding box pairs $P^p = \{(j,k)\}$, where $j \in I_b^p$ and $k \in I_b^q$. These intersecting bounding box pairs are used as the contact candidates for vesicle pairs in STIV calculation below.

*3.3.2. Phase 2: STIV computation.* We compute the STIV between the candidate vesicle pairs found by Alg. 1 to identify pairs that actually intersect.

Algorithm 2 summarizes the steps for computing the STIV. The algorithm starts by defining global index sets $I^p$ for vesicle points stored on each process (this requires an MPI reduction on the set of processes).

---

**Algorithm 2:** CONTACTVOLUME.

    **input** : Minimum separation distance $d_m$, vesicles' initial position $X_0^p$, candidate position $X_1^p$
    **output:** Contact Volumes $V^p$ and the Jacobian $J^p$

1  Compute $I^p$                                                                                      // MPIScan
2  $B^p \leftarrow \texttt{getBoundingBox}(X_0^p, X_1^p, d_m)$
3  $P^p \leftarrow \texttt{getIntersectingBBPair}(B^p, I^p)$                             // Alg. 1

4  $\overline{I}^p \leftarrow \texttt{getGhostVesicleID}(P^p)$
5  Send and receive $\{\overline{X}_0^p, \overline{X}_1^p\}$                                    // MPIAlltoallSparse
6  $\{V^p, I_v^p, J^p, \overline{J}^p\} \leftarrow \texttt{checkContact}(X_0^p, X_1^p, \overline{X}_0^p, \overline{X}_1^p, d_m)$

7  $\texttt{updateContactData}(I_v^p, J^p, \overline{J}^p)$                // MPIScan and MPIAlltoallSparse

---

Since vesicles are distributed over multiple processes, vesicles may have contact with vesicles on other processes; to compute STIV for all vesicles owned by a process $p$, positions and velocities of points of vesicles in $P^p$ need to be communicated to $p$. We refer to these copies of vesicle information as *ghost vesicles*.

Using the contact candidate vesicle pairs $P^p$, we compute the ghost vesicles' global index sets $\overline{I}^p$ which are used to distribute ghost vesicle point data $\{\overline{X}_0^p, \overline{X}_1^p\}$ (lines 4 and 5). On each process, using the contact detection method Section 3.1.2, we compute the contact volumes $V^p$, the contact volumes' index set $I_v^p$, the contact volume Jacobian $J^p$ with respect to the local vesicle points on process $p$, and the contact volume Jacobian $\overline{J}^p$ with respect to the ghost vesicle points copied to this process (line 6).

Initially, the contact volume index set $I_v^p$ is computed from local data and contains local indices. $I_v^p$ is converted to global indices by doing an MPI scan communication on $N_v^p$.

In addition, process $p$ sends back $\overline{J}^p$ to the owner process of each ghost vesicle and $J^p$ on each process is updated to store rows $J(I_v^p,:)$ of the global Jacobian matrix $J$.

In the contact volume calculation, each discretization point on the vesicle can only be involved in one contact volume, i.e., each column of $J$ has only a single nonzero element. As a consequence, $J^p$ can be compactly stored as vectors $g_X^p$ and $g_I^p$ where $g_X^p$ is the contact volume gradient with respect to $X_1^p$ and $g_I^p$ stores the contact volume indices the contact volume gradient components belong to. Next, we can proceed to solve the LCP and compute the magnitude of the contact force.

*3.3.3. Solving the linear complementarity problem.* The LCP matrix $\mathbf{B} = \mathbf{J}\mathbf{A}^{-1}\mathbf{G}\mathbf{J}^T$ is an $N_v \times N_v$ matrix, where $N_v$ is the number of contact volumes, $N_v = \mathcal{O}(N)$. Each entry $\mathbf{B}_{j,k}$ is the change in the $j^{\text{th}}$ contact volume induced by the $k^{\text{th}}$ contact force. Due to sparsity in matrices $\mathbf{J}$ and $\mathbf{A}$ (when the locally implicit scheme is used), the matrix $\mathbf{B}$ is sparse and typically diagonally dominant, since most STIV volumes are spatially separate. Two key algorithms are the parallel construction of the LCP matrix $\mathbf{B}$ and applying it to a vector (Alg. 3).

In Alg. 3, we form the LCP matrix $\mathbf{B}$, taking advantage of its sparsity and the sparsity of $\mathbf{J}$. After contact computation (Alg. 2), each process owns a list of contact volumes with indices $I_v^p$. A process $p$ stores rows $\mathbf{B}(I_v^p,:)$ of $\mathbf{B}$ and $\lambda(I_v^p)$, the local vector of Lagrange multipliers. Note that $I_v^p$ can be empty, which means there is no contact volume on process $p$. To form $\mathbf{B}(I_v^p,:)$, we loop over all the vesicles residing on process $p$.

For each contact volume pair $(j,k)$ that a vesicle $i$ is involved in, let $g_j$ and $g_k$ denote the rows $j$ and $k$ of the Jacobian matrix's columns with respect to vesicle $i$ sample point positions, then $\mathbf{B}(j,k)$ is updated as in line 7 of Alg. 3.

Finally, process $p$ needs to send $\mathbf{B}(j,k)$ to process $q$ if $j \in I_v^q$ and $p \neq q$, the owner of $j^{\text{th}}$ contact volume.

Since the matrix $\mathbf{B}$ and the vector $\lambda$ is distributed non-contiguously, some communication is needed to compute the matrix-vector product between $\mathbf{B}$ and $\lambda$. Process $p$ needs to send $\lambda(k)$ to process $q$ if there

---

**Algorithm 3:** FORM LCP MATRIX.

    **input** : Contact volume Jacobian $\mathbf{J}$ and CLI scheme matrix $\mathbf{A}$
    **output:** LCP matrix $\mathbf{B}$
    `// For` $p^{\text{th}}$ `process, form LCP matrix block` $\mathbf{B}(I_v^p,:)$`.  Use map data structure for sparse matrix B`

1  $\mathbf{B} \leftarrow \emptyset$
2  **for** *each local vesicle i* **do**
3     **for** *each contact volume j vesicle i involved in* **do**
4         **for** *each contact volume k vesicle i involved in* **do**
            `// accumulate change from` $k^{\text{th}}$ `contact force to` $j^{\text{th}}$ `contact volume`
5             $\mathbf{B}(j,k) = \mathbf{B}(j,k) + g_j \cdot \mathbf{A}(i,i)^{-1}\mathbf{G}(i,i)g_k$

6  **if** $j \notin I_v^p$ **then**
7     Send $\mathbf{B}(j,k)$ to process $q$ such that $j \in I_v^q$

is any non-zero entry $\mathbf{B}(j,k)$ where $j \in I_v^q$ and $k \in I_v^p$. After the communication, we can compute the LCP matrix-vector product $\mathbf{B}\lambda$ locally.

To solve the LCP, we use the minimum-map Newton method [5, Section 5.8], which only requires the application of the LCP matrix. Since the LCP matrix is already distributed, parallelizing the minimum-map Newton method LCP solver is straightforward. We presented a detail sequential minimum-map Newton method in [28].

## 4. Results

In this section, we present results characterizing the accuracy, robustness, and efficiency of a locally-implicit time stepping scheme (CLI) combined with our contact resolution framework in comparison to other schemes described in Section 3.2: with no contact resolution (i.e., LI scheme) and globally semi-implicit schemes with/without repulsion force (i.e., GI and RGI schemes).

- First, to demonstrate the robustness of our scheme in maintaining the prescribed minimum separation distance with different viscosity contrast $\nu$, we consider two vesicles in an extensional flow, Section 4.1.

- In Section 4.2, we explore the effect of minimum separation $d_m$ and its effect on collision displacement in shear flow. We demonstrate that the collision scheme has a minimum effect on the shear displacement.

- We present the timing for strong scalability and the weak scalability of our scheme.

- We close this section by reporting the computation cost (wall clock time) for simulations with different volume fractions, and compare the cost with RGI scheme.

Our experiments support the general observation that when vesicles become close, the LI scheme cannot, at a reasonable resolution of discretization, compute the interaction forces between close vesicles [49] and the time stepping becomes unstable. The GI scheme stays stable longer, but the iterative solver requires more and more iterations to reach the desired tolerance, which in turn implies higher computational cost for each time step. Finally, the RGI scheme requires choosing a penalty coefficient, which is typically done on case-by-case basis. If the penalty coefficient is too small, collisions many not be resolved, and an excessively large coefficient increases the error.

### 4.1. Extensional flow

Recall that to maintain accuracy of integral computations at a fixed surface resolution, one needs to ensure that a minimum separation distance is maintained. In order to demonstrate the ability of our framework to maintain a prescribed separation distance, we place two vesicles(with reduced volume 0.85) symmetrically with respect to the $z$ axis in the extensional flow $\boldsymbol{u} = [-x, y/2, z/2]$. For the experiments in this test, we use CLI and RGI schemes with first-order time stepping. We report both the minimal distance between two vesicles over time and the resulting minimal distance as a function of viscosity contrasts. Snapshots of the vesicle configuration in the CLI scheme are shown in Fig. 5.

In Fig. 6(a), we plot the distance between two vesicles over time using two different schemes (CLI and RGI). In RGI scheme, the vesicles get closer initially and fluctuate when the repulsion force is present. With different choices of repulsion force, RGI scheme may result in contact or large separation distance. On the other hand, our CLI scheme consistently maintains the desired minimum separation distance. In Fig. 6(b), we show the minimum distance between vesicles at the end of simulations, $T = 10$, as a function of viscosity contrast with no collision handling. We use adaptive time stepping GI scheme to run the simulations to get the final ($T = 10$) minimum distance between two vesicles. As expected, we observe the similar result as in 2D that the minimum distance between two vesicles decreases as the viscosity contrast is increased. In 3D, however, the minimum distance decreases much faster as the viscosity contrast is increased relative to 2D simulations [28].

To validate our estimates for the error due to piecewise-linear triangular approximation in the minimal separation calculation instead of the exact high-order geometry discretization, we plot the minimum distance at each step for two cases in Fig. 6(a): (i) the piecewise-linear triangular approximation; and (ii) the corresponding 4× upsampled shape. We observe that, as expected, the actual minimal distance for the smooth, high-order surface is smaller than the minimal distance for piecewise-linear triangular approximation, while the difference between two distances is small compared to the target minimum separation distance.

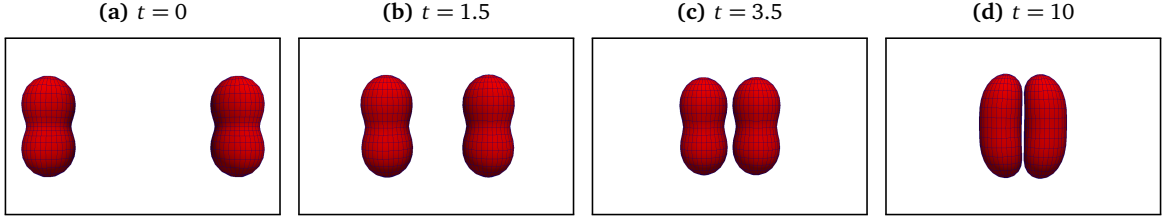

**Figure 5:** SNAPSHOTS OF TWO VESICLES IN EXTENSIONAL FLOW USING THE CLI SCHEME. *In the simulation, the viscosity contrast is set to 64. The CLI scheme maintains the desired minimum separation distance $d_m = 0.009$ as the distance between two vesicles decreases. With CLI scheme, two vesicles also maintain a symmetric configuration.*

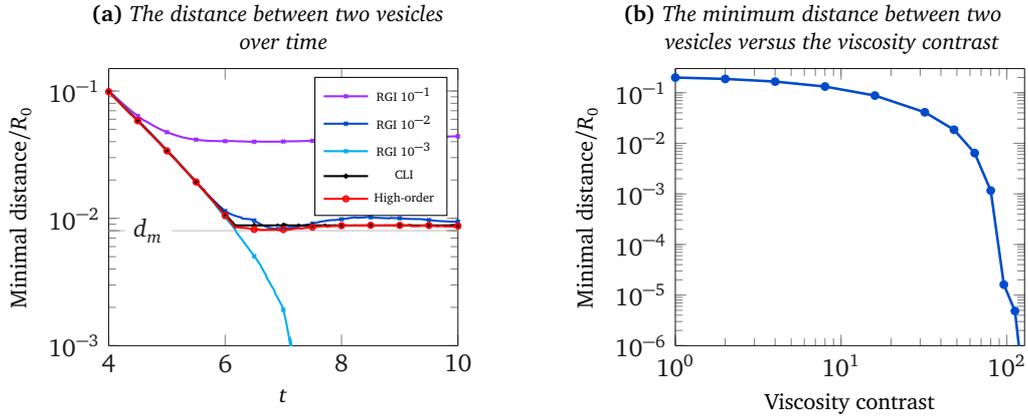

**Figure 6:** DISTANCE BETWEEN TWO VESICLES IN EXTENSIONAL FLOW. *(a) The minimal distance between two vesicles over time for CLI and RGI schemes with different choices of repulsion forces(in all case, the viscosity contrast is 64). The vesicles start to be in contact around $t = 6$. The black curve shows the minimal distance for the piecewise-linear triangular approximations to a $p = 32$ spherical harmonics grid. The red curve shows the estimated minimal distance between high order spectral surfaces (computed on a linear-triangulation to 4× upsampled surfaces). For the CLI scheme, the minimum separation distance $d_m$ is set to 0.009; as discussed in Section 3.1.2, we set separation distance to $(1 + 2\alpha)d_m$ with $\alpha = 0.05$ for $p = 32$ piecewise-linear triangular approximation for contact detection. For the RGI scheme, the repulsion coefficient $C_r$ is set to 0.1, 0.01 and 0.001. $R_0 := \sqrt{\text{Area}/4\pi}$ denotes the effective radius of a vesicle. In these simulations $R_0 = 1.136$ and it is used to normalize the minimal distance (the y-axis). The CLI scheme easily maintains the prescribed minimum separation of $d_m$. In this example, the RGI scheme maintains a similar minimum separation distance with $C_r = 0.01$, causes large separation distance with $C_r = 0.1$; and results in collision with $C_r = 0.001$. (b) With GI scheme, the final minimal distance between two vesicles as a function of viscosity contrast.*

With the minimum-separation constraint, any desired minimum separation distance between vesicles is maintained and the simulation is more robust as shown in Figs. 5 and 6. Moreover, the CLI scheme maintains the prescribed minimum separation, while the RGI scheme may fluctuate or collide depending on the prescribed repulsion coefficient $C_r$. We will show in Section 4.5 that with the prescribed repulsion coefficient of $C_r = 0.01$ and similar accuracy compared to the CLI scheme, the RGI does not prevent collisions and is more expensive than the CLI scheme.

### 4.2. Shear flow

To explore the effects of minimum separation distance on shear diffusivity, we place two vesicles of reduced volume 0.85 (to minimize the effect of vesicles' relative orientation on the dynamics) in an unbounded shear flow. We report the difference between centroids as a function of the minimum separation distance $d_m$ to demonstrate the effect of the minimum separation constrained system on the dynamics. For this experiment, we set the viscosity contrast to 5. Snapshot of the flow is shown in Fig. 7. In Fig. 8, we show the convergence of the vertical displacement between vesicle as a function of the minimum separation $d_m$ and the convergence rate of the scheme.

We consider two vesicles placed in a shear flow with (non-dimensional) shear rate $\chi = 1$. We consider a single shear rate in the experiments, with the observation that $\Delta t^{\text{stable}} \propto \chi^{-1}$ [48, Table 6] and [49, Table 4], the stable time step for other shear rates can be approximated from this. We use $\delta_t(d_m) := |z^1(t; d_m) - z^0(t; d_m)|$ to denote the vertical difference between two centroids of vesicles at time $t$. As shown in Fig. 7, at $T = 0$, two vesicles are placed with a vertical offset $\delta_0 = 0$. The background shear flow is $u = [\chi z, 0, 0]$ and the viscosity contrast is set to 5 in this experiment. In Fig. 8(a), we report the vertical offset between centroids over time as we increase the minimum separation distance. As the minimum separation parameter $d_m$ is decreased, the simulations with minimum-separation constraint converges to the precisely computed trajectory.

In Fig. 8(b), we report the convergence rate for the final error in centroid locations with some fixed minimum separation distance as we decrease the time step size. We use two $p = 16$ and $p = 32$ and $d_m$ accordingly. We observe first order convergence with CLI scheme. The LI scheme requires very small and often impractical time steps to prevent instability or intersection; we will revisit this in a later section.

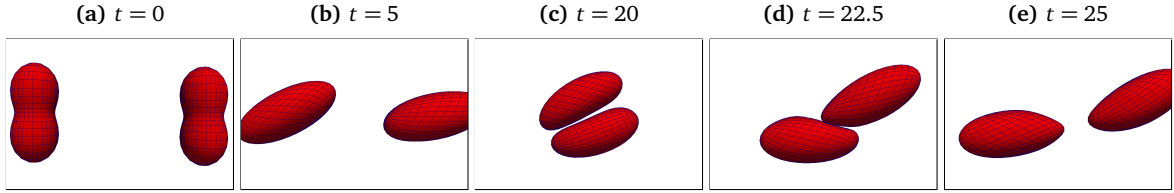

**Figure 7:** SHEAR FLOW EXPERIMENT. *The snapshots of two vesicles in shear flow. At $T = 0$, two vesicles are placed at $[-5.5, 0, 0]$ and $[0, 0, 0]$ respectively. The viscosity contrast for both vesicles is set to 5.*

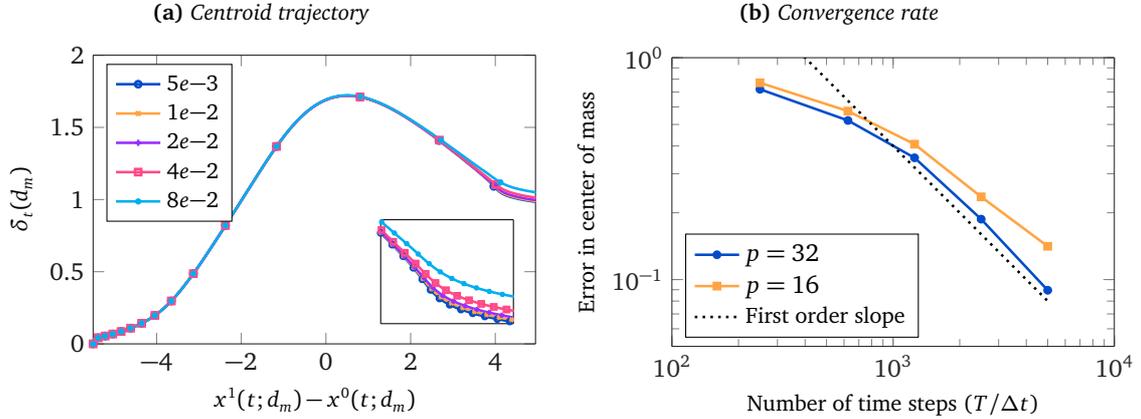

**Figure 8:** THE OFFSET $\delta_t(d_m)$ BETWEEN THE CENTROIDS OF TWO VESICLES IN SHEAR FLOW AND CONVERGENCE RATE. *The initial offset is $\delta_0 = 0$, viscosity contrast is set to 5, and the effective radius is $R_0 = 1.136$. (a) This plot shows the vertical offset $\delta_t(d_m) = |z^1(t) - z^0(t)|$ over time for different minimum separation distance $d_m$. The $d_m$ ranges from $5e{-}3$ to $8e{-}2$, the simulations converge, as we decrease the $d_m$. Spherical harmonic order $p = 16$ is used. (b) This plot shows the error in the final ($T = 25$) centroid location as we decrease the time step size for two spherical harmonic orders $p = 16$ and $p = 32$. We set $d_m = 5e{-}3$ for $p = 16$ and $d_m = 2.5e{-}3$ for $p = 32$. The final error in centroids is calculated with respect to the adaptive GI without repulsion force and $p = 32$ and error factor $E_f = 0.05$. As expected we observe first order convergence with our CLI scheme.*

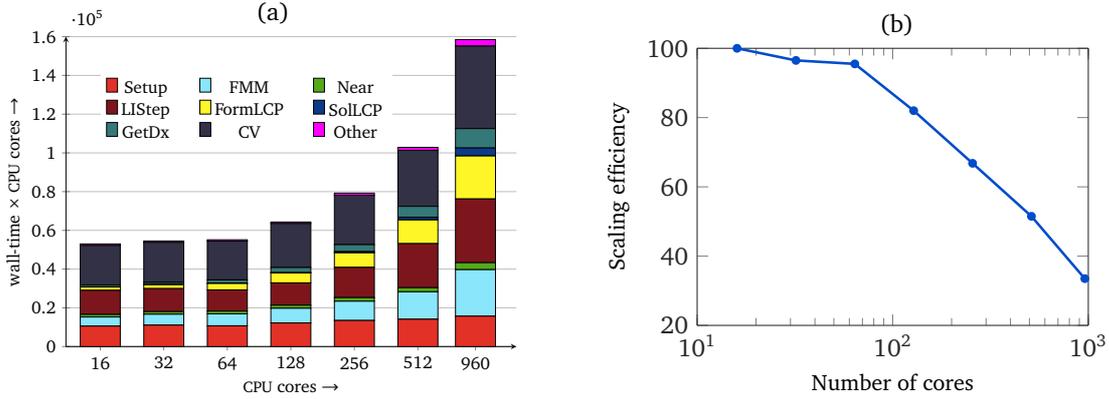

**Figure 9:** STRONG SCALABILITY FOR PERIODIC TAYLOR-VORTEX FLOW. *Setup is the setup phase for FMM, FMM is the actual FMM evaluation, Near is the near-singular calculation, LIStep represents the solve phase of LI scheme FormLCP and SolLCP are forming LCP matrix and solving LCP stages which constitute the time spent on LCP, CV is the contact volume calculation which corresponds to Alg. 2, GetDx is the phase of calculating new candidate solution, and Other represents all other calculations.*

*4.3. Strong scalability*

In this section and next, we use report the parallel scaling results for our framework. We used the Stampede1 system at the Texas Advanced Computing Center (TACC) to obtain the strong and weak scalability results. Each compute node in Stampede1 has two eight-core Intel Xeon E5-2680 CPUs running at 2.7GHz and 32GB of memory. We use the periodic Taylor vortex flow, as shown in Fig. 1, to investigate strong scaling. For this flow, the background velocity is

$$\boldsymbol{u}^\infty(x,y,z) = \alpha \sin\left(\frac{2\pi x}{L}\right)\cos\left(\frac{2\pi y}{L}\right)\sin\left(\frac{2\pi z}{L}\right)\boldsymbol{e}_1 + \alpha \cos\left(\frac{2\pi x}{L}\right)\sin\left(\frac{2\pi y}{L}\right)\sin\left(\frac{2\pi z}{L}\right)\boldsymbol{e}_2, \quad (4.1)$$

where $L$ is the periodic length and $\alpha$ is the scaling factor. For the simulations in this section, we choose $L = 28.4$ and $\alpha = 1$. Resulting timings are shown in Fig. 9.

We used the following simulation parameters:

- The number of vesicles is 1440; the vesicles are ellipsoidal, of effective radius $R_0 = 1.34$ with reduced volume 0.91 and the bending modulus is 0.1.

- Vesicle volume fraction is 58%.

- For spatial discretization, order $p = 16$ spherical harmonic were used, and the grid was upsampled to twice the resolution for collision detection.

- The time horizon is $T = 2$ and the time step size is $\Delta t = 0.1$.

- The block-diagonal solver relative tolerance is chosen to be $1e-5$.

The average number of contacts per vesicle stays about 2 per time step in our simulation. In Fig. 9, we report the total CPU time (wall-clock-time × CPU cores) for the number of cores ranging from 16 (1 compute node) to 960 (60 compute nodes). We achieve a speedup of 20.1 for the wall-clock-time or 33.5% strong scaling efficiency.

Figure 9 also shows the breakdown of the time spent in different parts of the code. As is evident in Fig. 9b, there are two main regimes in the scaling. When there are few vesicles per core (going from 64 to 128 cores), the load imbalance for the collision becomes more pronounced. The fraction of time spent on contact volume computation (CV) grows due to load imbalance with respect to the number of collisions and STIV calculation. Simple rebalancing by re-assigning different numbers of vesicles to processors is ineffective because of the subsequent increase in the load imbalance for the solver and FMM portions. The computational load of collision computations is dynamic and depends on the flow regime. To improve scalability of the collision computation in strong scaling regime, careful dynamic rebalancing, that takes into account the trade off between collision and solver balancing, will be needed.

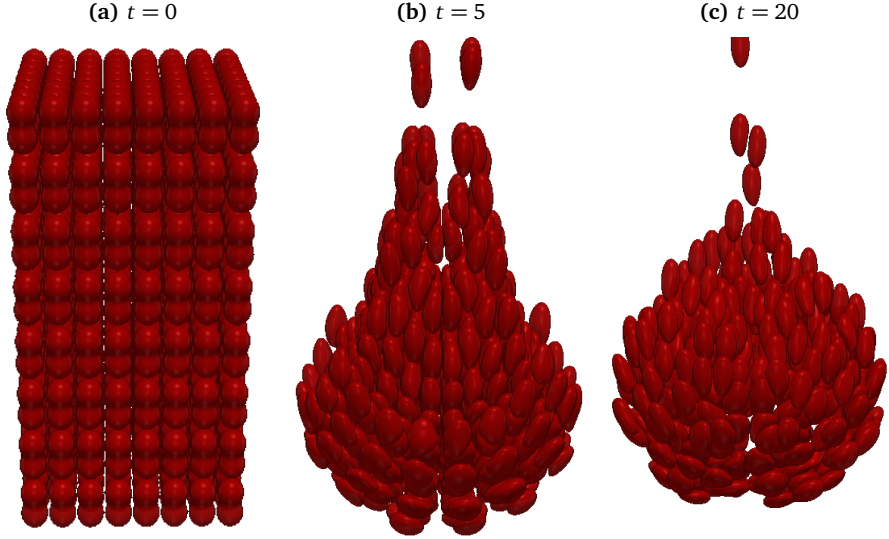

**Figure 10:** SEDIMENTATION OF POLY-DISPERSE VESICLES. *The snapshots of 512 (8 × 8 × 8 lattice) vesicles sedimenting under gravity. Each vesicle has a reduced volume of 0.91, bending modulus is in the range [0.05, 0.1], viscosity contrast in the range [0.5, 5] and an excess density of 1. We use time step size $\Delta t = 0.01$ and time horizon $T = 30$ in this simulation, the initial lattice has volume fraction of 53%.*

In our experiments, we observe that when the grain size (i.e., the number of vesicles per processor) is not too small, randomly distributing vesicles among processors with equal number of vesicles per processor achieves a good load balance for both solver and collision parts of the code.

The solver time (LIStep) dominates for very small grain size, since the linear system is block-diagonal and each process requires a different number of GMRES iterations for convergence; processes with fewer GMRES iterations will wait for processes requiring more GMRES iterations.

*4.4. Weak Scalability*

To showcase different flows, we use the sedimentation of a poly-disperse suspension of vesicles as in Fig. 10 on 16K CPU cores for our weak scaling study. The scaling results are shown in Fig. 11. We use time step size of $\Delta t = 0.01$ and a time horizon $T = 0.1$. All other simulation parameters are equal to those of the strong scaling test.

We present two sets of results for 1 vesicle per core (Fig. 11 left) and 8 vesicles per core (Fig. 11 right). We present a breakdown of the time spent on different functions of our algorithm as we scale from 16 cores to 16K cores. Similar to the strong scaling case, load imbalance with respect to the number of collisions causes the timing to grow. In this flow, different regions of space have very different number of collisions.

Another factor that affects the timing is the number of contact-resolving iterations, which grows from 4 to 8 as we increase the number of cores from 16 to 16K (i.e., as the problem size increases). For uniform lattice as in sedimentation experiment Fig. 10, the number of contact-resolving iterations stabilizes to about 8 as we increase the number of cores.

*4.5. High volume-fraction flows*

In our final experiment, we investigate the effectiveness of our scheme in modeling flows with high volume-fraction $\phi$. We use the wall-clock-time for a fixed time horizon to quantify the cost of simulation for each scheme (CLI, RGI, and time-adaptive RGI).

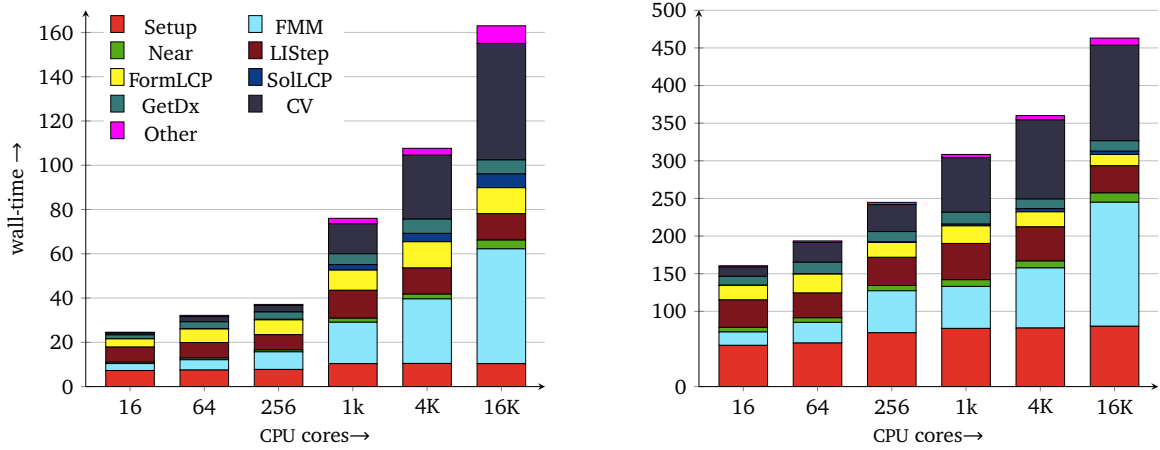

**Figure 11:** WEAK SCALABILITY RESULTS FOR POLY-DISPERSE SEDIMENTATION. *The left figure is with grain size* 1 *vesicle per core, the right figure is with grain size* 8 *vesicles per core. The flow snapshots are shown in Fig. 10 and the flow parameters are outlined therein. For each case, we present a breakdown of the wall-clock-time spent on each of the different functions in our algorithm. See the caption of Fig. 9 for the description of the labels.*

For this experiment, we use Taylor-vortex flow with 168 prolate vesicles distributed on a staggered lattice as shown in Fig. 12. The shape and distribution of vesicles are chosen to achieve high volume fractions), the simulation time horizon is set to $T = 15$. We change the periodic length $L$ and spacing between vesicles to obtain different volume fractions. All of the simulations are executed on a dedicated node with the same type of CPU with 1 MPI process to ensure the wall time used is calculated consistently.

We run three sets of experiments:

1. *CLI scheme:* We use the CLI scheme to run the simulations with different volume fractions $\phi$ and different time step sizes $\Delta t$. The minimum separation $d_m$ is set to 0.009. We report the total wall-clock-time (in seconds) in Table 4a. Since the time stepping scheme is locally implicit, large time steps cause the simulation to diverge. We mark those cases by "LI-div". As expected, the wall-clock-time shows weak dependence on the volume fraction since there are more collisions.

2. *Adaptive time stepping with repulsion:* We use adaptive time stepping RGI scheme to run simulations with different error factors $E_f$, different volume fractions $\phi$, and a fixed repulsion coefficient $C_r = 0.01$, which maintains similar separation distance as $d_m = 0.009$ for the CLI scheme above, Fig. 6(a). In the adaptive scheme, the error factor is the tolerance for the error committed in each simulation time unit [30]. If the problem is non-stiff, one expects the time step to be proportional to $E_f$. Therefore, the counterpart of $E_f$ in non-adaptive schemes is the step size $\Delta t$. In Table 4b, we report the wall-clock time in seconds provided that the simulation finishes. If $\Delta t$ of the adaptive scheme is reduced below ($1e-11$) we abort the simulation (these are marked as "Limit" in the table). If a collision occurs, we stop the simulation and mark it as "Col".

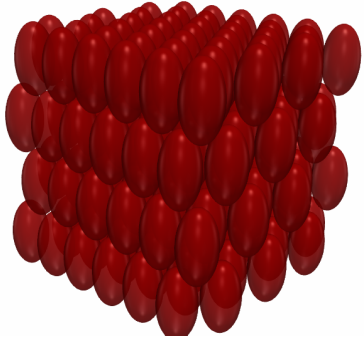

**Figure 12:** VOLUME FRACTION EXPERIMENT. *An example of the initial distribution of 168 vesicles in Taylor-vortex flow. We modify the spacing between vesicles to obtain different volume fractions.*

In Table 4b, we chose the widest possible range for $E_f$ (going from .25 to 4.0) to present the full picture with respect to the simulation cost. Comparing the results in Table 4, we see that the adaptive scheme is one to two orders of magnitude more expensive for similar cases. For example, for $\phi = 0.55$, the RGI scheme with $E_f = 0.5$ requires $143e3$ seconds compared to $2.7e3$ of the CLI scheme with $\Delta t = 0.4$; a 53× speedup.

3. *Non-adaptive RGI scheme:* To compare the cost of RGI scheme with that of CLI, we report the wall-clock time for fixed $\Delta t$ and different repulsion coefficients $C_r$; the volume fraction for this experiment is fixed at $\phi = 0.5$. The time steps are chosen to include those used for the CLI scheme in Table 4. In Table 5 we report the wall-clock time (in seconds) when the simulation finishes. The cases where the scheme diverges, because the chosen time step is unstable due to the stiffness introduced by the repulsion force, are marked as "RGI-div".

Larger repulsion coefficients avoid collision at the expense of accuracy. In [30, Figure 5(f)], the error analysis for different repulsion coefficients shows that large repulsion coefficient will introduce significant error in the center of mass trajectory. The largest repulsion coefficient we test here is 0.10 which already has a significant error. To match the error of CLI, $C_r$ needs to be set around 0.01 that in turn requires very small time step.

The results in Table 5 show that for a fixed $\Delta t$, the repulsion coefficient needs to be adjusted for the simulation to succeed. There are many factors influencing the choice of repulsion coefficient, e.g., bending modulus, viscosity contrast, volume fraction, reduced volume, background flow and vesicle shape, which makes an automatic choice difficult.

Comparing to the row corresponding to volume fraction $\phi = 0.5$ of the CLI scheme in Table 4, the CLI scheme is at least 5× faster than RGI with the "right" choice of the repulsion coefficient.

| $\phi$ | $\Delta t$ | | | |
|---|---|---|---|---|
| | 0.1 | 0.2 | 0.4 | 0.8 |
| 0.35 | $3.1e3$ | $2.1e3$ | $1.8e3$ | LI-div |
| 0.40 | $3.3e3$ | $2.2e3$ | $1.9e3$ | LI-div |
| 0.45 | $3.7e3$ | $2.5e3$ | $2.1e3$ | LI-div |
| 0.50 | $4.2e3$ | $2.7e3$ | $2.3e3$ | LI-div |
| 0.55 | $5.2e3$ | $3.4e3$ | $2.7e3$ | LI-div |

**(a)** *CLI scheme*

| $\phi$ | $E_f$ | | | | |
|---|---|---|---|---|---|
| | 0.25 | 0.5 | 1.0 | 2.0 | 4.0 |
| 0.35 | $150e3$ | $150e3$ | $142e3$ | $136e3$ | Col |
| 0.40 | $139e3$ | $135e3$ | $133e3$ | $127e3$ | Col |
| 0.45 | Limit | $172e3$ | $173e3$ | $173e3$ | Col |
| 0.50 | Limit | $158e3$ | $158e3$ | $154e3$ | Col |
| 0.55 | Limit | $143e3$ | $139e3$ | Col | Col |

**(b)** *Adaptive RGI scheme*

**Table 4:** WALL-CLOCK-TIME VS. VOLUME FRACTION. *Wall-clock-time (seconds) for simulating 168 vesicles in Taylor vortex flow with different volume fractions. (a) For the CLI scheme, $d_m = 0.009$. (b) For the RGI scheme, the $C_r = 0.01$.*

| $C_r$ | $\Delta t$ | | | | | |
|---|---|---|---|---|---|---|
| | 0.025 | 0.05 | 0.1 | 0.2 | 0.4 | 0.8 |
| 0.01 | $49.7e3$ | Col | Col | Col | Col | RGI-div |
| 0.02 | $46.6e3$ | $33.3e3$ | $20.2e3$ | Col | Col | RGI-Div |
| 0.05 | $42.7e3$ | $30.1e3$ | $18.8e3$ | $11.1e3$ | Col | RGI-Div |
| 0.10 | $40.1e3$ | $29.3e3$ | $17.5e3$ | $10.5e3$ | Col | RGI-Div |

**Table 5:** WALL-CLOCK-TIME VS. REPULSION COEFFICIENT. *Wall-clock-time (seconds) for the non-adaptive RGI scheme with different repulsion coefficient $C_r$ and fixed time step $\Delta t$. The volume fraction is fixed at $\phi = 0.5$.*

## 5. Conclusion

We have introduced new parallel algorithms for efficient 3D simulation of non-dilute suspensions of deformable particles immersed in Stokesian fluid. We demonstrated the parallel scaling of the algorithms on up to 16K CPU cores. Moreover, we demonstrated through numerical experiments that our scheme is orders of magnitude faster than the alternatives for several setups.

We would like to thank George Biros, David Harmon, Dhairya Malhotra, Matthew Morse, Bryan Quaife, Michael Shelley, and Etienne Vouga for stimulating conversations about various aspects of this work. This work was supported by the US National Science Foundation (NSF) through grants DMS-1320621 and DMS-1436591.

## Appendix A. Space-time volume

The space-time volume is a 3D volume embedded in 4D space parameterized by $\phi$, $\theta$ and $t$, that is, by the map

$$P(\phi, \theta, t) = (x(\phi, \theta, t), y(\phi, \theta, t), z(\phi, \theta, t), \epsilon t) = (X, \epsilon t),$$

where $\epsilon$ is a scaling factor and has the unit of velocity. It is used to make the units consistent. Let $g_{\alpha\beta}$ and $h_{\alpha\beta}$ respectively denote the 3D and 4D metric tensors, i.e.,

$$g_{\alpha\beta} = X_\alpha \cdot X_\beta, \quad \alpha, \beta \in \{\phi, \theta\},$$
$$h_{\alpha\beta} = P_\alpha \cdot P_\beta, \quad \alpha, \beta \in \{\phi, \theta, t\},$$

and $g = \det(g_{\alpha\beta})$, $h = \det(h_{\alpha\beta})$. Elements $h_{\alpha\beta}$ are related to $g_{\alpha\beta}$, spatial velocity, and basis of tangent space as

$$h = \det(h_{\alpha\beta}) = \det \begin{bmatrix} g_{\phi\phi} & g_{\phi\theta} & u \cdot X_\phi \\ g_{\theta\phi} & g_{\theta\theta} & u \cdot X_\theta \\ u \cdot X_\phi & u \cdot X_\theta & \epsilon^2 + u \cdot u \end{bmatrix}.$$

Expanding and simplifying the expression for the determinant, we have

$$h = g\left[\epsilon^2 + (u \cdot u) - (u \cdot X^\alpha)(u \cdot X_\alpha)\right],$$

where $X^\alpha := g^{\alpha\beta} X_\beta$ is the contravariant basis and $g^{\alpha\beta}$ denotes the dual tensor to $g_{\alpha\beta}$. By the definition of the contravariant basis, we can decompose the velocity as

$$u = (u \cdot X^\alpha) X_\alpha + (u \cdot n) n.$$

This in turn implies,

$$u \cdot u = (u \cdot X^\alpha)(u \cdot X^\beta) X_\alpha \cdot X_\beta + (u \cdot n)^2 = (u \cdot X^\alpha)(u \cdot X_\alpha) + (u \cdot n)^2.$$

Combining these results, we have

$$h = g\left[\epsilon^2 + (u \cdot n)^2\right].$$

The space-time volume is now easily computed by

$$V^C = \int \sqrt{h} \, d\phi \, d\theta \, dt = \int \sqrt{g(\epsilon^2 + (u \cdot n)^2)} \, d\phi \, d\theta \, dt = \int \sqrt{\epsilon^2 + (u \cdot n)^2} \, dA \, dt,$$

where we substituted $dA = \sqrt{g} \, d\theta \, d\phi$.

## Appendix B. Boundary integral formulation

We follow the standard potential theory [38] and express the solution of Eq. (2.6) as a system of singular integro-differential equations on all vesicle surfaces. The solution of Eq. (2.6) is expressed by the combination of single- and double-layer integrals. We denote the single-layer integral on the vesicle surface $\gamma_i$ by

$$\mathcal{G}_{\gamma_i}[f](x) := \int_{\gamma_i} \mathbf{G}(x - Y) \cdot f(Y) \, dA, \tag{B.1}$$

where $f$ is an appropriately defined density. The double-layer integral is defined as

$$\mathcal{T}_{\gamma_i}[q](x) := \int_{\gamma_i} n(Y) \cdot \mathbf{T}(x - Y) \cdot q(Y) \, dA, \tag{B.2}$$

where $q$ is the appropriately defined density.

In 3D, the fundamental solutions of the Stokes equation are given by the Stokeslet tensor $\mathbf{G}$ and the Stresslet tensor $\mathbf{T}$,

$$\mathbf{G}(r) = \frac{1}{8\pi\mu} \frac{1}{\|r\|} \left( I + \frac{r \otimes r}{\|r\|^2} \right), \qquad \mathbf{T}(r) = -\frac{3}{4\pi\mu} \frac{r \otimes r \otimes r}{\|r\|^5}.$$

For any point $x \in \mathbb{R}^3$, the velocity can be expressed as the superposition of velocities due to each vesicle

$$\alpha u(x) = u^\infty(x) + \sum_i u^i(x), \quad x \in \mathbb{R}^3, \quad \alpha = \begin{cases} 1 & x \in \mathbb{R}^3 \setminus \gamma, \\ \nu_i & x \in \varpi_i, \\ (1 + \nu_i)/2 & x \in \gamma_i, \end{cases} \tag{B.3}$$

where $u^\infty(x)$ represents the background velocity field, $u^i$ denotes the velocity contributions from vesicle $i$, and $\nu_i = \mu_i/\mu$ denotes the viscosity contrast of the $i^{\text{th}}$ vesicle. We define the complementary velocity for each vesicle as $\bar{u}^i = \alpha u - u^i$.

The velocity induced by the $i^{\text{th}}$ vesicle is expressed as an combination of single- and double-laryer integrals[41]:

$$u^i(x) = \mathcal{G}_{\gamma_i}[f](x) + (1 - \nu_i)\mathcal{T}_{\gamma_i}[u](x) \qquad (x \in \mathbb{R}^3), \tag{B.4}$$

where the double-layer density $u$ is the total velocity from Eq. (B.3) and $f$ is the traction jump across the vesicle membrane. The traction jump on each vesicle is equal to the sum of bending, tensile, and collision forces Eq. (2.7)

$$f(X) = f_b + f_\sigma + f_c \qquad (X \in \gamma). \tag{B.5}$$

Substitute Eq. (B.4) into Eq. (B.3), we get the equation for the interfacial velocity of each vesicle,

$$\frac{(1 + \nu_i)}{2} u(X) = \bar{u}^i(X) + \mathcal{G}_{\gamma_i}[f](X) + (1 - \nu_i)\mathcal{T}_{\gamma_i}[u](X) \qquad (X \in \gamma_i), \tag{B.6}$$

subject to local inextensibility

$$\nabla_\gamma \cdot u = 0 \qquad (X \in \gamma). \tag{B.7}$$